
\magnification1200
\footline={\hss{\vbox to 2cm{\vfil\hbox{\rm\folio}}}\hss}
\nopagenumbers
  \def\rightheadline{{\hfil{\sevenrm
  On the multiplicity of zeros of the zeta-function}\hfil\tenrm\folio}}

  \def\leftheadline{{\tenrm\folio\hfil{\sevenrm
   Aleksandar Ivi\'c }\hfil}}
  \def\emptyheadline{\hfil}
  \headline{\ifnum\pageno=1 \emptyheadline\else
  \ifodd\pageno \rightheadline \else \leftheadline\fi\fi}

\font\teneufm=eufm10
\font\seveneufm=eufm7
\font\fiveeufm=eufm5
\newfam\eufmfam
\textfont\eufmfam=\teneufm
\scriptfont\eufmfam=\seveneufm
\scriptscriptfont\eufmfam=\fiveeufm
\def\mathfrak#1{{\fam\eufmfam\relax#1}}

\font\tenmsb=msbm10
\font\sevenmsb=msbm7
\font\fivemsb=msbm5
\newfam\msbfam
      \textfont\msbfam=\tenmsb
      \scriptfont\msbfam=\sevenmsb
      \scriptscriptfont\msbfam=\fivemsb
\def\Bbb#1{{\fam\msbfam #1}}

\def \CC {\Bbb C}

\def \RR {\Bbb R}

\def\DJ{\leavevmode\setbox0=\hbox{D}\kern0pt\rlap
 {\kern.04em\raise.188\ht0\hbox{-}}D}
\def\b{\beta} \def\g{\gamma} \def\R{\Re{\rm e}\,}
\def\z{\zeta} \def\r{\rho} \def\d{{\,\rm d}} \def\G{\Gamma}
\def\s{\sigma}
\def\hf{{\textstyle{1\over2}}}
\def\txt#1{{\textstyle{#1}}}
\font\cc=cmss10 scaled 1200
\font\dd=cmcsc10 scaled 1200
\font\ee=cmr9

\centerline{\sevenbf Acad\'emie Serbe des  Sciences et des Arts, Bulletin
CXVIII Classe des Sciences }

\smallskip\centerline{\sevenbf
Math\'ematiques et naturelles, Sciences Math\'ematiques No. 24, Belgrade 1999,
119--132}

\bigskip
\centerline{\cc ON THE MULTIPLICITY OF ZEROS OF THE ZETA-FUNCTION}

\bigskip\bigskip\centerline{\dd  A. IVI\'C }

\bigskip\bigskip
\centerline {(Presented at the 1st Meeting, held on February  26, 1999)}
\bigskip\bigskip
{\ee A b s t r a c t. Several results are obtained concerning
multiplicities of zeros of the
Riemann zeta-function $\zeta(s)$. They include upper bounds for
multiplicities, showing that zeros with large multiplicities have
to lie to the left of the line $\s = 1$. A zero-density counting function
involving multiplicities is also discussed.


}
\bigskip\bigskip
\centerline{\dd 1. Introduction}
\bigskip\bigskip

Let $r = m(\rho)$ denote the multiplicity of the complex
zero $\rho = \b + i\g$
of the Riemann zeta-function $\z(s)$. This means that $\z(\r) = \z'(\r)
= \ldots = \z^{(r-1)}(\r) = 0$, but $\z^{(r)}(\r) \not = 0$. All known zeros
$\r$ are simple (i.e. $m(\r) = 1$), and it may well be that they are all
simple, although the proof of this is
certainly beyond reach at present. In estimating $m(\r)$ one may suppose
that $\b \ge \hf$ and that $\g > 0$, since $1-\r$ and $\overline\r$
 are zeros of $\z(s)$ if $\r$ is a zero.

\smallskip
It  seems that there exist no good upper bounds
in the literature for $m(\r)$. All that
appears to be known unconditionally is
$$
m(\b + i\g) \;\ll\; \log \g.\eqno(1.1)
$$
On the Lindel\"of Hypothesis (LH) that $\z({1\over2} +it) \ll
 |t|^\varepsilon$
this bound can be improved to
$$
m(\b + i\g) \;=\; o(\log \g)\qquad(\g \to \infty),\eqno(1.2)
$$
and on the Riemann Hypothesis (RH) that $\r = {1\over2} + i\g$ to
$$
m(\b + i\g) \;\ll\; {\log \g\over\log\log\g}.\eqno(1.3)
$$
Furthermore, on the RH H.L. Montgomery [4] proved
that at least 2/3 of the zeros are simple.
It transpires that the estimation of $m(\b + i\g)$ is a very difficult
problem, and one which is not satisfactorily solved even under the
assumption of the LH or the RH.
To see how one obtains (1.1)--(1.3) recall that
for $N(T)$, the number of zeros $\b+ i\g$ for which $0 < \g \le T$, one
has the classical Riemann-von Mangoldt formula (see [1] and [6])
$$
N(T) = {T\over2\pi}\log\left({T\over2\pi}\right) - {T\over2\pi} + {7\over8}
+ S(T) + O\left({1\over T}\right),\quad S(T) = {1\over\pi}\arg\z({1\over2}
+ iT),\eqno(1.4)
$$
where $\arg\z({1\over2}+ iT)$ is obtained by continuous
variation along the straight lines
joining the points $2, 2 + iT, {1\over2} + iT$, starting with the value 0.
If $T$ is the ordinate of a zero lying on the critical line,
 then $S(T) =S(T + 0)$. One has (see [6]) the bounds
$$
S(T) \ll \log T,\quad S(T) = o(\log T) \quad({\rm LH}),\quad
S(T) \ll {\log T\over\log\log T}\quad ({\rm RH}),
$$
and these bounds combined with the trivial inequality
$$
m(\b + i \g) \;\le\; N(\g + H) - N(\g - H)  \qquad(0 < H \le 1)\eqno(1.5)
$$
easily yield (1.1)--(1.3), respectively.
It seems, however, that these estimates
are much too large, and that perhaps one even has
$$
m(\b + i\g) \;\ll_\varepsilon \;(\log\log\g)^{1+\varepsilon},\eqno(1.6)
$$
which is of course still much weaker than the conjecture that all zeros
are simple.
The use of pointwise estimates for $S(T)$ certainly cannot give anything
close to (1.6), since one has
$$
S(T) = \Omega_\pm\left(\left({\log T\over\log\log T}\right)^{1/3}\right),
\quad
S(T) = \Omega_\pm\left(\left({\log T\over\log\log T}\right)^{1/2}\right)\quad
({\rm RH}),
$$
proved by K.-M. Tsang [7] and H.L. Montgomery [5], respectively.
(As usual, $f = \Omega_\pm(g)$ means that $\limsup_{x\to\infty}f(x)/g(x)
= +\infty$ and $\liminf_{x\to\infty}f(x)/g(x) = -\infty$ both hold).
 One could use (1.5) with $H = o(1)\;(\g\to\infty)$ to
try to improve (1.1)--(1.3). In view of (1.4) this is equivalent to
obtaining bounds for $S(\g + H) - S(\g - H)$, but no satisfactory
results seem to be known for this problem.

\smallskip
In this note we shall seek other approaches to the estimation of
$m(\b + i \g) $, and several bounds will be proved in the sequel.
We shall also discuss a zero-density counting function involving
multiplicities.
 In what follows $C$ will denote positive, absolute constants,
not necessarily the same ones at each occurrence.

\bigskip\bigskip
\centerline{\dd 2. Bounds for multiplicities}
\bigskip\bigskip

In this section we shall formulate and prove the  results involving
upper bounds for  $m(\b + i\g)$. We start with

\bigskip
THEOREM 1. {\it If $\z(\b + i\g) = 0, \hf < \b < 1$ and $\g \ge \g_0 > 0$,
 then}
$$
m(\b+i\g) \le {1\over\log{1\over2-2\b}}\left(\max_{\s\ge{1\over2},|t|
\le{1\over2}}\log|\z(\s+i\g+it)| + O(\log\log\g)\right).
\eqno(2.1)
$$

{\bf Proof of Theorem 1.}  We shall use Jensen's classical formula
(see e.g. [2], pp. 257-258). Namely if $f(z)$ is regular in $|z| \le R$
 and $f(0) \not = 0$, then
$$
\sum_{f(\rho)=0,|\rho|\le R}\log{R\over|\rho|} + \log|f(0)| =
{1\over2\pi}\int_0^{2\pi}\log|f(Re^{i\theta})|\d\theta,\eqno(2.2)
$$
where the zeros $\r$ are counted according to their multiplicities. In (2.2)
we take $f(z) = \z(1+i\g+z), R = \hf, \g \ge \g_0, \hf <\b < 1$. Then
$f(0) \not = 0$, in fact (see Lemma 12.3 of [1])
$$
\z(\s + it)  \gg \log^{-2/3}t(\log\log t)^{-1/3}
\quad(\s \ge 1 - C(\log t)^{-2/3}(\log\log t)^{-1/3}, t  \ge t_0 > 0).
\eqno(2.3)
$$
If $r = m(\b+i\g)$, then (2.2) gives
$$\eqalign{
r\log{\hf\over1-\b} + \log|\z(1+i\g)| &\le {1\over2\pi}\int_0^{2\pi}
\log|\z(1+i\g+\hf e^{i\theta})|\d \theta\cr&
\le \max_{\s\ge{1\over2},|t|\le{1\over2}}\log |\z(\s+i\g+it)|,}
$$
and using (2.3) we immediately obtain (2.1). A slight improvement is possible
with the choice $f(z) = \z(1 - \delta(\g) +z), R = \hf - \delta(\g),\,
\delta(\g) = C(\log \g)^{-2/3}(\log\log\g)^{-1/3}$.

\bigskip
THEOREM 2. {\it If $\z(\b + i\g) = 0, \hf < \b < 1, \g \ge \g_0 > 0$
and $\,c\,$ is a constant satisfying $c > 1 - \b$, then}
$$
m(\b+i\g) \;\le\; {{c+\b-1\over c+\b-\hf}\,\max\limits_{|t|\le\log^2\g}
\log|\z(\hf+i\g+it)| + O(\log\log\g)  \over\log\left\{
{c\over1-\b}\left({\b-\hf\over c}\right)^{c+\b-1\over c+\b-
{1\over2}}\right\}}.
\eqno(2.4)
$$

{\bf Proof of Theorem 2.} Let $r = m(\b+i\g)$ and $\cal D$ be the rectangle
with vertices $\hf - \b\pm i\log^2\g,\,c\pm i\log^2\g$.
If $X\,(0 < X \ll \g^C)$ is a parameter which will be suitably chosen,
then by the residue theorem we obtain
$$
{\z(1-\b+\r)\over(1-\b)^r} \;=\; {1\over2\pi i}
\int_{\cal D}X^{s-1+\b}\G(s-1+\b){\z(s+\r)\over s^r}\d s
\quad(\rho = \b+i\g).\eqno(2.5)
$$
By using (2.3) and Stirling's formula for the gamma-function it follows
from (2.5) that
$$\eqalign{
{1\over(1-\b)^r\log\g} &\ll e^{-\log^2\g} + X^{-{1\over2}}
\max_{|t|\le\log^2\g}\,{|\z(\hf+i\g+it)|\over(\b-\hf)^r}\cr&
+ X^{c-1+\b}\int_{-\log^2\g}^{\log^2\g}|\G(c-1+\b+it)|\,
{|\z(c+\b+i\g+it)|\over|c+it|^r}\d t\cr&
\ll (\b-\hf)^{-r}X^{-{1\over2}}M + X^{c-1+\b}c^{-r} + e^{-\log^2\g},}
\eqno(2.6)
$$
where
$$
M \;:=\; \max_{|t|\le\log^2\g}\,|\z(\hf+i\g+it)|.
$$
We choose
$$
X \;=\; M^{1\over c+\b-{1\over2}}\left({c\over\b-\hf}\right)
^{r\over c+\b-{1\over2}}
$$
to equalize the terms containing $X$ in (2.6), noting that in view of
(1.1) $\,X \ll \g^C\,$ will hold.
Then after taking logarithms we obtain (2.4) from (2.6).

\bigskip
To assess the strength of (2.4) and compare it with (2.1) of Theorem 1,
take $c = {\txt{3\over2}} - \b$. Then (2.4) gives
$$
m(\b+i\g) \;\le\;{\log M\over\log\left\{({3\over2}-\b)(\b-{1\over2})
\over(1-\b)^2\right\}} + O_\b(\log\log\g).\eqno(2.7)
$$
But we have
$$
{1\over\log\left\{({3\over2}-\b)(\b-{1\over2})
\over(1-\b)^2\right\}} \;\le\;{1\over\log{1\over2-2\b}}
$$
for $\b \ge {5-\sqrt{5}\over4} = 0.69058\ldots \,$, hence in this range
essentially (2.7) improves (2.1), and the range can be further increased by an
appropriate choice of the constant $c$.

\bigskip
Both Theorem 1 and Theorem 2 could be modified to cover the case
$\b = \hf$, i.e. the zeros on the critical line. In that case the maxima
in (2.1) and (2.4) would have to include zeta-values on the lines left
of $\s =\hf$, which would not yield good results in view of the
functional equation
$$
\z(s) = \chi(s)\z(1 -s),\qquad \chi(s) \asymp t^{{1\over2}-\s}
\quad(s = \s + it,\,t \ge t_0 > 0).
$$
The main merit of  Theorem 1 and Theorem 2 is that they
show that $m(\b+i\g)$ is small if
$\b$ is close to the line $\s=1$.
Namely, if $\b \ge 1 - \log^{-\delta}\g
\;(0 < \delta \le {2\over3}; \delta > {2\over3}$ is not possible
in view of (2.3)),
we obtain either from Theorem 1 or Theorem 2 that
$$
m(\b + i\g) \;\ll_\delta\; {\log\g\over\log\log\g},
$$
which is the same as (1.3). Therefore zeros close to the line $\s=1$,
 if they exist,
must have small multiplicities. On the other hand, if there exist zeros
with large multiplicities (``large" in the sense that $m(\b + i\g)/
\log\log\g\to\infty$), then they must be far away from the line $\s=1$.
This is precisely given by

\medskip
THEOREM 3. {\it If $\,\z(\b+ i\g) = 0$, then for }$\g \ge \g_0 > 0$
$$
m(\b + i\g) \;\ll\; (1 - \b)^{3/2}\log\g + \log\log\g,\eqno(2.8)
$$
{\it and if $\,\lim\limits_{\g\to\infty}{m(\b + i\g)\over
\log\log\g} = +\infty$, then there is a constant $C > 0$
such that for} $\g \ge \g_0 > 0$
$$
\b \;\le\;1 \;-\;C\left({m(\b + i\g)\over\log\g}\right)^{2/3}.\eqno(2.9)
$$

\bigskip

\bigskip {\bf Proof of Theorem 3.} Let $\cal D$ be the rectangle
with vertices $-2(1-\b)\pm i\log^2\g,\,1\pm i\log^2\g$. We can
suppose that $\b \ge {\txt{3\over4}}$, for otherwise (2.8) is trivial in
view of (1.1), and so is (2.9) if $C$ is sufficiently small.
 For $\b \ge {\txt{3\over4}}$
formula (2.5) (with the above $\cal D$)
is valid, and to estimate the integrand on the line $\s = -2(1-\b)$
we shall use the bound
$$
\G(w)  \;\ll\; {e^{-|{\rm Im}\,w|}\over |w|}.
$$
Since $(1-\b)^{-1} \ll \log\g$ it follows from (2.5) that
$$
{(1-\b)^{-r}\over\log\g} \ll e^{-\log^2\g} + 2^{-r}(1-\b)^{-r}
X^{-3(1-\b)}\log\g\,\max_{|\g|\le\log^2\g}\,|\z(3\b-2+i\g+it)| + X^\b.
\eqno(2.10)
$$
To bound the zeta-factor in (2.10) we use the estimate (see Ch. 6 of [1])
$$
\z(\s+it) \;\ll\; t^{C(1-\s)^{3/2}}\log^{2/3}t\qquad(t \ge t_0,\,
\hf \le \s \le 1,\, C \le 122),
$$
which is the strongest known one when $\s$ is sufficiently close to
$\s = 1$.
It follows that
$$
{1\over(1-\b)^r} \ll \left(2^{-r}(1-\b)^{-r}X^{-3(1-\b)}\g^{C_1(1-\b)^{3/2}}
+ X^\b\right)\log^{C_2}\g\qquad(C_1,\,C_2 > 0).
$$
To equalize the terms containing $X$ we choose
$$
X \;=\; 2^{-{r\over3-2\b}}(1-\b)^{-{r\over3-2\b}}\g^{C_1(1-\b)^{3/2}
\over3-2\b}.
$$
The condition $X\ll \g^C$ will hold in view of (2.1), and we obtain
$$
\left\{2^{\b\over3-2\b}(1-\b)^{3\b-3\over3-2\b}\right\}^r
\;\ll\;\g^{C(1-\b)^{3/2}}\log^{C_2}\g.
$$
Taking logarithms we have
$$
r(\b\log2 +(3\b-3)\log(1-\b)) \;\le\; C(1-\b)^{3/2}\log\g + O(\log\log\g).
\eqno(2.11)
$$
Now since $\b \ge {\txt{3\over4}} $ and
$$
(3\b-3)\log(1-\b) \;\ge \;0\qquad(0 < \b < 1),
$$
we obtain (2.8) from (2.11), and if additionally
$$
\lim_{\g\to\infty}\,{r\over\log\log\g} \;=\;
\lim_{\g\to\infty}\,{m(\b+i\g)\over\log\log\g} = +\infty,\eqno(2.12)
$$
then (2.9) also follows from (2.11) (or from (2.8)).
Thus if (2.12) holds for some zero $\r$, then
(2.9) shows that $\r$ lies to the left of the sharpest known
zero-free region for $\z(s)$ implied by (2.3).

\medskip
The condition (2.12) appears significant in another context. Namely
if one assumes that the zeros near the line $\s = 1$ are isolated from
one another, then one can improve the known zero-free region of
$\z(s)$ (see Ch. 6 of [1]) that
$$
\b \;\le \; 1 - C(\log\g)^{-2/3}(\log\log\g)^{-1/3}
\qquad(\z(\b+i\g) = 0,\,\g \ge \g_0 > 0).
$$
This was done by N. Levinson [3], who proved the following result:
{\it If for some $\delta > 0$ and sufficiently large constant $T_0 > 0$
the zeros $\rho$ which lie in $\b > 1 - \delta,\,|\g| > T_0$ are all
isolated in the sense that there is no zero of $\z(s)$ other than
$\rho = \b + i\g$ in the rectangle ($s = \s + it$)
$$
1 - \delta \;<\;\s\;<\;1,\quad |t - \g| \;<\;2\delta,\eqno(2.13)
$$
then there are no zeros of $\z(s)$ in the region
$$
\s \;>\;1 \;-\;{C\over\log\log t}\qquad(\,|t| > T_0\,)\eqno(2.14)
$$
for suitable $C > 0$. Although isolated, these zeros need not be simple.}

\medskip
By going through Levinson's proof and making the appropriate modifications
one can obtain an upper bound for $\b$ depending on $m(\b+i\g)$, and
which shows that under he condition (2.12) the zero-free region (2.14)
can be improved. This is

\medskip
THEOREM 4. {\it Suppose that for some $\delta > 0$ and
sufficiently large constant $T_0 > 0$
the zeros $\rho$ which lie in $\b > 1 - \delta,\,|\g| > T_0$ are all
isolated in the sense that there is no zero of $\z(s)$ other than
$\rho = \b + i\g$ in the rectangle} (2.13). {\it If $r = m(\b + i\g)$
for such a zero, then}
$$
\b \;\le\;1 \;-\; {r^{1\over2\log\log\g}\over C\log\log\g}.\eqno(2.15)
$$

\smallskip
{\bf Proof of Theorem 4.} The theorem actually shows that $r = m(\b + i\g)$
cannot be larger than $(C\log\log\g)^{\log\log\g}$,
for otherwise (2.15) would yield that $\b < 1-\delta$, which
is impossible. One can rewrite (2.15) as
$$
 m(\b + i\g) \;\le\;\left(C(1 - \b)\log\log\g\right)^{2\log\log\g},
$$
which is the upper bound for $m(\b+i\g)$ under the ``isolation hypothesis".
To prove the theorem, let ($\Lambda(n)$ is the von Mangoldt function)
$$
F_m(s) := (-1)^m{{\rm d}^{m-1}\over {\rm d}s^{m-1}}\,
\left({\z'(s)\over\z(s)}\right)
= \sum_{n=1}^\infty\Lambda(n)n^{-s}\log^{m-1}n\quad(\s =
{\rm Re\,}s > 1;\,m \ge 1).
$$
We use (see (1.43) of [1]; here exceptionally in the next two
formulas $\g$ is Euler's constant and
not the imaginary part of a zeta-zero)
$$
{\z'(s)\over\z(s)} = \log(2\pi) - 1 - \hf\g - \hf{\G'({s\over2}+1)
\over\G({s\over2}+1)} + \sum_\rho\left({1\over s - \rho} +
{1\over\rho}\right),
$$
$$
{\G'(z)\over \G(z)} \;=\; -\g + \sum_{n=0}^\infty\left({1\over n+1}
- {1\over n+z}\right),
$$
(2.13) and $N(T + 1) - N(T) \ll \log T$ to obtain
$$
F_m(s) = -\sum_{\rho;|\g-t|\le 1}(s - \rho)^{-m} + O(\log|t|) \quad(s =
\s + it \in \CC,\,m \ge 2),\eqno(2.16)
$$
where the $O$--constants throughout the proof are uniform in $m$.
Therefore the isolation hypothesis  yields
$$
F_m(\s + i\g) = - {r\over(\s-\b)^m} + O\left({\log\g\over\delta^m}\right)
\qquad(\,r = m(\b+i\g);\,m \ge 2).\eqno(2.17)
$$
In the rectangle $1 - \delta < \s < 1,|t - 2\g| < \delta$, there is at
most one zero of $\z(s)$ by the isolation hypothesis. Consider the
case when there is such a zero, say $\b_1 + i\g_1$; if there is no such
zero (2.20) will hold with $r_1 = 0$. Then in place of (2.17) one
obtains from (2.16)
$$
F_m(\s + 2i\g) = - {r_1\over(\s-\b_1+i(2\g-\g_1))^m} +
O\left({\log\g\over\delta^m}\right) \qquad(\,r_1 = m(\b_1+i\g_1);\,m \ge 2)
\eqno(2.18)
$$
and, for $\s > 1$, one also has
$$
F_m(\s) \;=\; {1\over(\s-1)^m} + O(1)\qquad(\,m \ge 2\,).\eqno(2.19)
$$
From the definition of $F_m(s)$ and the classical inequality
$$
3 + 4\cos\theta + \cos2\theta = 2(1 + \cos\theta)^2 \ge 0
\qquad(\theta \in \RR)
$$
it follows, for $\s > 1$, that
$$
\R\left\{3F_m(\s) + 4F_m(\s + it) + F_m(\s + 2it)\right\} \ge 0,
$$
which in view of (2.17)--(2.19) yields, for $\s > 1$ and $m \ge 2$,
$$
{3\over(\s-1)^m} - {4r\over(\s-\b)^m} - \R
{r_1\over(\s_1-\b_1+i(2\g-\g_1))^m} \ge -{C\log\g\over\delta^m}.\eqno(2.20)
$$
Now one chooses
$$
M = [\log\log\g] + 1,\quad \s = 1 + 100M(1-\b),
$$
and the chief feature of Levinson's method is that $m$ in (2.20) can be
chosen in such a way that $m \asymp M$ and that ${\rm Re\,}(\cdots)$ is
non-negative. As shown in detail in [3], there are two cases: in the
first case $m = M \,(\ge 2)\,$ suffices for the proof, and in the
second case $M \le m < 4M$ will hold. Since the analysis is identical
with the one needed in our proof, it will be omitted here. With the
above choices of $M$ and $\s$ one obtains from (2.20) in the case
$m = M$ (the other case is quite analogous)
$$
(100M(1-\b))^{-M}\left\{3 - 4r\left({100M\over1+100M}\right)^M\right\}
\ge -{C\log\g\over\delta^M}.\eqno(2.21)
$$
But we have
$$\eqalign{
3 - 4r\left({100M\over1+100M}\right)^M &= 3 - 4r\left(1 + {1\over100M}
\right)^{-M} < 3 - 4re^{-{1\over100}}\cr&
< 3 - {7r\over2} \le -{r\over2}.\cr}
$$
Consequently (2.21) yields
$$
{r\delta^M\over2C\log\g} \;\le \; (100M(1-\b))^M,
$$
which easily leads to (2.15). The choice of $M$ is made so that
$(\log\g)^{1/M}$ is bounded.

\bigskip
We also remark  that
there is a possibility to bound $m(\b+i\g)$, provided one has a good
lower bound of the form
$$
\int_\delta^{2\delta}|\z(\b+i\g+i\alpha)|^k\d \alpha \,\ge\, \ell\,
= \,\ell(\g,\delta,k)\qquad(0 < \delta < {\txt{1\over4}},
\;\g \ge \g_0 > 0) \eqno(2.22)
$$
for $k = 1,2$. Namely for $\b \ge \hf$ let $\cal D$ be the
rectangle with vertices
${\textstyle 1\over4} -\b\pm i\log^2\g, 2 \pm i\log^2\g, \z(\rho) = 0,
\rho = \b+i\g,
\g \ge \g_0 > 0$, and let $0 < \alpha \le 1$. Then by the residue  theorem
$$
{\z(\b+i\g+i\alpha)\over(i\alpha)^r}
= {1\over2\pi i}\int_{\cal D} \Gamma(s-i\alpha){\z(s+\rho)\over s^r}\d s.
$$
This gives
$$
\z(\b+i\g+i\alpha) \ll \alpha^r\left(\g(\b - {\textstyle {1\over4}})^{-r}
+ 2^{-r}\right) \ll \alpha^r\g(\b - {\textstyle {1\over4}})^{-r},
$$
and consequently, if $\delta$ is a constant satisfying $0 < \delta
< {1\over4}$, we have
$$
\int_\delta^{2\delta}|\z(\b+i\g+i\alpha)|^k\d \alpha
\ll \g^k(\b - {\textstyle {1\over4}})^{-rk}\int_\delta^{2\delta}\alpha^{rk}
\d\alpha \ll \delta^{rk}\g^k.
$$
Therefore by using (2.22) and taking logarithms we obtain
$$
m(\b+i\g) = r \le {1\over\log
\left({1\over\delta}\right)}
\left(\log\g - {1\over k}\log\ell + O(1)\right).\eqno(2.23)
$$
Hence (2.23) shows that the upper bound for $m(\b+i\g)$ can be made to
depend on $\ell$ in (2.22). We would like to let $\delta \to 0+$ in (2.23)
and obtain (1.2). However,
by using the argument on top of p. 219 of [6]
and the first  inequality on p. 230, it follows that in
(2.22) one can take $\ell = \delta \gamma^{-A/\delta}$
or even $\delta \gamma^{A\log\delta}$
for some absolute $A > 0$. These bounds, unfortunately, are
too weak to yield (1.2).

\bigskip\bigskip
\centerline{\dd 3. The zero-density counting function}
\bigskip\bigskip

\medskip
A problem related to the estimation of $m(\b+i\g)$
is to estimate $N^{(r)}(\s,T)$, the number
of zeros $\rho = \b + i\g$ of $\z(s)$ such that $m(\rho) \ge r,
\b \ge \s \,(\ge \hf)$ and $|\g| \le T$. Note that $N^{(1)}(\s,T)
\equiv N(\s,T)$, where as usual $N(\s,T)$ denotes the number of zeros
$\rho$ such that $\b \ge \s$ and $|\g| \le T$, counted with their
multiplicities. To deal with this problem we define,
for $x,c > 0$ and integral $r \ge 0$,
$$
f_r(x) \;:=\; {1\over2\pi i}\int_{(c)}\,\G(s)x^{-s}s^{-r}\d s,\eqno(3.1)
$$
where as usual $\int_{(c)} = \lim_{T\to\infty}\int_{c-iT}^{c+iT}$. We have
$f_0(x) = e^{-x}$ by the Mellin inversion formula of the gamma-integral, and
since the integral in (3.1) is absolutely convergent we obtain
$$
f_r'(x) = -{1\over2\pi i}\int_{(c)}\G(s)x^{-s-1}s^{1-r}\d s
= -x^{-1}f_{r-1}(x)\qquad(r \ge 1).\eqno(3.2)
$$
For $c \ge 1$ we have, by (3.1) and trivial estimation, $f_r(x) \ll_c x^{-c}$,
hence by integration of (3.2) we obtain
$$
f_r(x) = \int_x^\infty f_{r-1}(t){\d t\over t}\qquad(r \ge 1).\eqno(3.3)
$$
By induction on $r$ we see from (3.2) and (3.3) that $f_r(x)$ is positive
and monotonically decreasing in $(0,\,\infty)$ (and $\lim_{x\to0+}f_r(x)
= +\infty$ for $r \ge 1$). From (3.2) and (3.3) we obtain, using repeatedly
integration by parts,
$$
\leqalignno{f_r(x)& = f_{r-1}(t)\log{t\over x}\Bigg|_{x}^\infty -
\int_x^\infty f_{r-1}'(t)\log {t\over x}\d t\cr&
= \int_x^\infty f_{r-2}(t)\log {t\over x}{\d t\over t} =
{1\over2!}\int_x^\infty f_{r-3}(t)\log^2 {t\over x}{\d t\over t} = \ldots
&(3.4)\cr&
= {1\over r!}\int_x^\infty \log^r({t\over x})e^{-t}\d t\qquad(r \ge 1).\cr}
$$
From (3.3) we obtain, by induction on $r$,
$$
f_r(x) \;\le\;x^{-r}e^{-x}\quad(x > 0,\,r =  0,1,2,\ldots\,).\eqno(3.5)
$$
From (3.4) we have
$$
f_r(x) \ge {1\over r!}\int_1^\infty \log^r\bigl({t\over x}\bigr)e^{-t}\d t
\ge {\log^r({1\over x})\over r!}\int_1^\infty e^{-t}\d t =
{\log^r({1\over x})\over r!e}\quad(0 < x \le 1, r \ge 1).\eqno(3.6)
$$
Set ($\mu(n)$ is the M\"obius function)
$$
M_X(s) = \sum_{n\le X}\mu(n)n^{-s}\quad(X \ge 2),\qquad a(n)
= \sum_{d|n,d\le X}\,\mu(d),
$$
and note that $a(1) = 1, a(n) = 0$ for $2 \le n \le X$ and
$|a(n)| \le d(n)$, the number of divisors of $n$.
Let  $r = R + 1 =m(\r) = m(\b + i\g) \to \infty, \b > \hf,
\g \ge \g_0 > 0$, and suppose that
for some constant $B > 0$
$$
1 \;\ll\;X\;\le\;Y,\quad Y^{1-\b} \;\ll\; \g^B.\eqno(3.7)
$$
Then by (3.1) and absolute convergence we obtain, since $\G(s + 1)  = s\G(s)$,
$$\eqalign{
\sum_{n=1}^\infty a(n)n^{-\r}f_R({n\over Y})
&= \sum_{n=1}^\infty a(n)n^{-\r}\cdot{1\over2\pi i}\int_{(2)}
\Bigl({Y\over n}\Bigr)^s\G(s)s^{-R}\d s\cr&
= {1\over2\pi i}\int_{(2)}\z(\r + s)M_X(\r + s)Y^s\G(s + 1)s^{-r}\d s.
\cr}\eqno(3.8)
$$
Note that by (3.5) we have, as $Y \to \infty$, uniformly in $R$
$$
\eqalign{&\sum_{n=1}^\infty a(n)n^{-\r}f_R({n\over Y}) =
a(1)f_R({1\over Y}) + \sum_{n>X}a(n)n^{-\r}f_R({n\over Y})\cr&
= f_R({1\over Y}) + \sum_{X<n\le 2Y\log Y}a(n)n^{-\r}f_R({n\over Y})
+ O\left(\sum_{n>2Y\log Y}d(n)n^{-\b}e^{-n/Y}\right)\cr&
= f_R({1\over Y}) + \sum_{X<n\le 2Y\log Y}a(n)n^{-\r}f_R({n\over Y})
+ o(1).\cr}\eqno(3.9)
$$
Now if $\rho$ is counted by $N^{(r)}(\s,T)$, then
(3.6), (3.8) and (3.9) give that either $|\g| \le \log^2T$ or
$$
1 \;\ll\; \Big|\sum_{X<n\le2Y\log Y}b(n)n^{-\rho}\Big|\qquad
\left(b(n) = a(n){f_R({n\over Y})\over f_R({1\over Y})}\right)
$$
or
$$
1 \;\ll\; {R!\over\log^RY}\int_{-\infty}^\infty
|\z(\hf + it + i\g)M_X(\hf + it + i\g)|e^{-\pi|t|/2}
Y^{1/2-\b}\,{\d t\over|{1\over2} - \b + it|^r}.
$$
From the properties of $f_r(x)$ it follows that $|b(n)| \le a(n)$, and
the above bounds are analogous to (11.9) and (11.10) of [1], which appear
as the starting point in the derivation of upper bounds for $N(\s,T)$.
The  difference is that, if $r \ge 2$ is fixed, then in the analogue of
(11.10) of [1] we shall have an additional factor $\log^{-R}Y = \log^{1-r}Y$.
Hence eventually for $N^{(r)}(\s,T)$ we shall obtain the same upper bound as
we would for $N(\s,T)$, only it will be smaller by a factor of $\log^AT,\,
A = A(r,\s) > 0$ (for $r \ge 2$ fixed). In each specific upper bound
for $N^{(r)}(\s,T)$ this constant $A = A(r,\s)$ can be explicitly
evaluated. Actually one expects that, for $\hf \le \s < 1$ fixed,
$$
\lim_{T\to\infty}\,{N^{(r)}(\s,T)\over N(\s,T)} \;=\;0\qquad(r \ge 2)
\eqno(3.10)$$
will hold, although proving this is out of reach at present. Note that the
function $f_r(x)$, defined by (3.1), is not the only kernel function which
may be used to estimate $N^{(r)}(\s,T)$, but the use of other similar
functions does not appear to yield sharper results. One could also use the
above method to obtain a pointwise estimate for $m(\b+i\g)$
(e.g., by choosing $X = 2Y\log Y$ in (3.8)), but the bound
that will be obtained will not be better than (2.1).

\medskip
One may compare the conjectural (3.10) with the known bound (see p. 246
of [6])
$$
N_r(T) \;\ll\; N(T)e^{-C\sqrt{r}},\eqno(3.11)
$$
where $N(T)$ denotes the number of zeros $\rho$ with $0 < \g \le T$
and $N_r(T)$ the number  of zeros $\rho$ with $0 < \g \le T$
and $m(\b+i\g) = r$. The methods used in obtaining (3.11) involve
moments of the function $S(t + h) - S(t)$ and do not seem to be able
to yield (3.10).

\bigskip
\centerline{\cc REFERENCES}
\bigskip

\item{[1]} A. Ivi\'c, {\it The Riemann zeta-function}, John Wiley
and Sons, New York, 1985.

\smallskip

\item{[2]} S. Lang, {\it Complex Analysis (2nd ed.)}, Springer, GTM 103,
New York etc., 1985.

\smallskip
\item{[3]} N. Levinson, {\it Zeros of the Riemann zeta-function
near the 1-line}, J. Math. Anal. Appl. {\bf 25}(1969), 250-253.

\smallskip
\item{[4]} H.L. Montgomery, {\it The pair correlation of zeros of the
zeta-function}, Proc. Symp. Pure Math. {\bf 24}, AMS, Providence R.I.,
1973, 181-193.

\smallskip
\item{[5]} H.L. Montgomery, {\it Extreme values of the Riemann
zeta-function}, Comment. Math. Helv. {\bf 52}(1977), 511-518.

\smallskip\item{[6]} E.C. Titchmarsh, {\it The theory of the Riemann
zeta-function (2nd ed.)},  Oxford, Clarendon Press, 1986.

\smallskip
\item{[7]} K.-M. Tsang, {\it Some $\Omega$--theorems for the Riemann
zeta-function}, Acta Arithmetica {\bf 46}(1986), 369-395.

\bigskip\bigskip\bigskip\bigskip
\noindent \sevenrm Aleksandar Ivi\'c, Katedra Matematike

\noindent  RGF-a Universiteta u Beogradu, \DJ u\v sina 7

\noindent 11000 Beograd, Serbia and Montenegro

\noindent e-mail:
\sevenbf aivic@matf.bg.ac.yu, \enskip ivic@rgf.bg.ac.yu
\bye